\documentclass{article}

\title{\LARGE \textbf{On Dirac's Conjecture}}
\author{Zh.G. Nikoghosyan}

\begin{document}

\maketitle

\begin{abstract}
Let $G$ be a 2-connected graph, $l$ be the length of a longest path in $G$ and $c$ be the circumference - the length of a longest cycle in $G$. In 1952, Dirac proved that $c>\sqrt{2l}$ and conjectured that $c\ge 2\sqrt{l}$. In this paper we present more general sharp bounds in terms of $l$ and the length $m$ of a vine on a longest path in $G$ including Dirac's conjecture as a corollary: if $c=m+y+2$ (generally, $c\ge m+y+2$) for some integer $y\ge 0$, then $c\ge\sqrt{4l+(y+1)^2}$ if $m$ is odd; and $c\ge\sqrt{4l+(y+1)^2-1}$ if $m$ is even. \\

Key words: longest cycle, longest path, circumference.

\end{abstract}

\section{Introduction}

We consider only undirected graphs with no loops or multiple edges. Let $G$ be a 2-connected graph. We use $c$ and $l$ to denote the circumference (the length of a longest cycle) and the length of a longest path of $G$. A good reference for any undefined terms is \cite{[1]}.

In 1952, Dirac \cite{[3]} proved the following.  \\

\noindent\textbf{Theorem A} \cite{[3]}. In every 2-connected graph, $c>\sqrt{2l}$. \\

In the same paper \cite{[3]}, Dirac conjectured a sharp version of Theorem A.\\

\noindent\textbf{Conjecture A} \cite{[3]}. In every 2-connected graph, $c\ge2\sqrt{l}$.\\

In this paper we present more general sharp bounds in terms of $l$ and the length of a vine on a longest path of $G$ including Dirac's conjecture as a corollary. In order to formulate this result, we need some additional definitions and notations.

The set of vertices of a graph $G$ is denoted by $V(G)$ and the set of edges by $E(G)$. If $Q$ is a path or a cycle, then the length of $Q$, denoted by $l(Q)$, is $|E(Q)|$ - the number of edges in $Q$. We write a cycle $Q$ with a given orientation by $\overrightarrow{Q}$. For $x,y\in V(Q)$, we denote by $x\overrightarrow{Q}y$ the subpath of $Q$ in the chosen direction from $x$ to $y$. We use $P=x\overrightarrow{P}y$ to denote a path with end vertices $x$ and $y$ in the direction from $x$ to $y$. We say that vertex $z_1$ precedes vertex $z_2$ on $\overrightarrow{Q}$ if $z_1$, $z_2$ occur on $\overrightarrow{Q}$ in this order, and indicate this relationship by $z_1\prec z_2$. We will write $z_1\preceq z_2$ when either $z_1=z_2$ or $z_1\prec z_2$.

Let $P=x\overrightarrow{P}y$ be a path. A vine of length $m$ on $P$ is a set
$$
\{L_i=x_i\overrightarrow{L}_iy_i: 1\le i\le m\}
$$
of internally-disjoint paths such that\\

$(a)$ $V(L_i)\cap V(P)=\{x_i,y_i\} \ \ (i=1,...,m)$,

$(b)$ $x=x_1\prec x_2\prec y_1\preceq x_3\prec y_2\preceq x_4\prec ...\preceq x_m\prec y_{m-1}\prec y_m=y$ on $P$.\\

\noindent\textbf{Theorem 1}. Let $G$ be a 2-connected graph. If $\{L_1,L_2,...,L_m\}$ is a vine on a longest path of $G$ and $c=m+y+2$ for some integer $y\ge 0$, then

$$
c\ge\left\{
\begin{array}{lll}
\sqrt{4l+(y+1)^2} & \mbox{when} &  m\equiv 1(mod\ 2), \\ \sqrt{4l+(y+1)^2-1} & \mbox{when} & %
m\equiv 0(mod\ 2).
\end{array}
\right.
$$

The following lemma guarantees the existence of at least one vine on a longest path in a 2-connected graph.\\

\noindent\textbf{The Vine Lemma} {[2]}. Let $G$ be a $k$-connected graph and $P$ a path in $G$. Then there are $k-1$ pairwise-disjoint vines on $P$.\\

\section{The proof of Theorem 1}

Let $P=x\overrightarrow{P}y$ be a longest path in $G$ and let 
$$
\{L_i=x_i\overrightarrow{L}_iy_i: 1\le i\le m\}
$$
be a vine of length $m$ on $P$. Put

$$
L_i=x_i\overrightarrow{L}_iy_i \ \ (i=1,...,m), \ \ A_1=x_1\overrightarrow{P}x_2, \ \ A_m=y_{m-1}\overrightarrow{P}y_m,
$$
$$
A_i=y_{i-1}\overrightarrow{P}x_{i+1} \ \ (i=2,3,...,m-1), \ \  B_i=x_{i+1}\overrightarrow{P}y_i  \  \  (i=1,...,m-1),
$$
$$
l(A_i)=a_i \  \ (i=1,...,m),\  \ l(B_i)=b_i \  \ (i=1,...,m-1).
$$

By combining appropriate $L_i, A_i, B_i$, we can form the following cycles:
$$
Q_0=\bigcup_{i=1}^m(A_i\cup L_i),
$$
$$
Q_j=\bigcup_{i=j+1}^{m-j}(A_i\cup L_i)\cup B_j \cup B_{m-j} \  \   (j=1,2,...,\lfloor (m-1)/2\rfloor).
$$

Since $l(L_i)\ge 1$ \ $(i=1,2,...,m)$ and $a_1\ge 1$, $a_m\ge 1$, we have
$$
c\ge l(Q_0)= \sum_{i=1}^m l(L_i)+a_1+a_m +\sum_{i=2}^{m-1}a_i\ge m+2,
$$
Let $c=m+y+2$, where $y\ge 0$. If $a_1+a_m\ge y+3-\sum_{i=2}^{m-1}a_i$, then
$$
l(Q_0)= \sum_{i=1}^m l(L_i)+(a_1+a_m)+\sum_{i=2}^{m-1}a_i
$$
$$
\ge m+\left(y+3-\sum_{i=2}^{m-1}a_i\right)+\sum_{i=2}^{m-1}a_i= m+y+3>c,
$$
a contradiction. Hence,
$$
a_1+a_m\le y+2-\sum_{i=2}^{m-1}a_i.    \eqno{(1)}
$$
Next, if $b_1+b_{m-1}\ge y+5-\sum_{i=2}^{m-1}a_i$, then
$$
l(Q_1)=\sum_{i=2}^{m-1} l(L_i)+(b_1+b_{m-1})+\sum_{i=2}^{m-1}a_i
$$
$$
\ge (m-2)+\left(y+5-\sum_{i=2}^{m-1}a_i\right)+\sum_{i=2}^{m-1}a_i=m+y+3>c,
$$
a contradiction. Hence,
$$
b_1+b_{m-1}\le y+4-\sum_{i=2}^{m-1}a_i\le y+4.
$$
Furthermore, if $b_2+b_{m-2}\ge y+7-\sum_{i=3}^{m-2}a_i$, then
$$
l(Q_2)=\sum_{i=3}^{m-2} l(L_i)+(b_2+b_{m-2})+\sum_{i=3}^{m-2}a_i
$$
$$
\ge (m-4)+\left(y+7-\sum_{i=3}^{m-2}a_i\right)+\sum_{i=3}^{m-2}a_i=m+y+3>c,
$$
again a contradiction. Hence
$$
b_2+b_{m-2}\le y+6-\sum_{i=3}^{m-2}a_i\le y+6.
$$
Analogously,
$$
b_j+b_{m-j}\le y+2(j+1)-\sum_{i=j+1}^{m-j}a_i
$$
$$
\le y+2(j+1) \ \  (j=1,2,...,\lfloor(m-1)/2\rfloor).        \eqno{(2)}
$$

\textbf{Case 1}. $m=2k+1$ for some integer $k\ge 0$.

By (2),
$$
\sum_{i=1}^{m-1}b_i=\sum_{j=1}^{(m-1)/2}(b_j+b_{m-j})\le \sum_{j=1}^{(m-1)/2}(y+2(j+1))
$$
$$
=\frac{m-1}{2}(y+2)+2\sum_{j=1}^(m-1)/2=\frac{m-1}{2}(y+2)+\frac{(m-1)(m+1)}{4}.
$$
Then
$$
l=(a_1+a_m)+\sum_{i=2}^{m-1}a_i+\sum_{i=1}^{m-1}b_i
$$
$$
\le \left(y+2-\sum_{i=2}^{m-1}a_i\right)+\sum_{i=2}^{m-1}a_i+ \frac{m-1}{2}(y+2)+\frac{(m-1)(m+1)}{4}
$$
$$
=(y+2)\frac{m+1}{2}+\frac{(m-1)(m+1)}{4}
$$
$$
= (y+2)\frac{c-y-1}{2}+\frac{(c-y-3)(c-y-1)}{4}=\frac{c^2-(y+1)^2}{4},
$$
implying that $c\ge \sqrt{4l+(y+1)^2}$.\\

\textbf{Case 2}. $m=2k$ for some integer $k\ge 1$.

As in Case 1,
$$
\sum_{i=1}^{m-1}b_i=\sum_{j=1}^{(m-2)/2}(b_j+b_{m-j})+b_{m/2}
$$
$$
\le(y+2)\frac{m-2}{2}+\frac{(m-2)m}{4}+b_{m/2}.
$$
For the cycle 
$$
Q^\ast=B_{m/2}\cup B_{(m-2)/2}\cup A_{m/2}\cup L_{m/2}
$$
we have
$$
l(Q^\ast)=b_{m/2}+b_{(m-2)/2}+a_{m/2}+l(L_{m/2})\le c=y+m+2,
$$
implying that
$$
b_{m/2}+b_{(m-2)/2}\le y+m+1.
$$
Then it easy to see that
$$
\sum_{i=1}^{m-1}b_i\le (y+2)\frac{m-2}{2}+\frac{(m-2)m}{4}+\frac{y}{2}+\frac{m+2}{2}.
$$

Furthermore, we have
$$
l=(a_1+a_m)+\sum_{i=2}^{m-1}a_i+\sum_{i=1}^{m-1}b_i
$$
$$
\le (y+2)+(y+2)\frac{m-2}{2}+\frac{(m-2)m}{4}+\frac{y}{2}+\frac{m+2}{2}
$$
$$
=\frac{(m+2)^2+2y(m+1)}{4}=\frac{(c-y)^2+2y(c-y-1)}{4}=\frac{c^2-y^2-2y}{4},
$$
implying that $c\ge \sqrt{4l+(y+1)^2-1}$. Theorem is proved.   \quad \quad \quad  \quad \rule{7pt}{6pt}\\

Let $P=x\overrightarrow{P}y$ be a path and let 
$$
\{L_i=x_i\overrightarrow{L}_iy_i: 1\le i\le m\}
$$
be a vine of length $m$ on $P$. Put

$$
L_i=x_i\overrightarrow{L}_iy_i \ \ (i=1,...,m), \ \ A_1=x_1\overrightarrow{P}x_2, \ \ A_m=y_{m-1}\overrightarrow{P}y_m,
$$
$$
A_i=y_{i-1}\overrightarrow{P}x_{i+1} \ \ (i=2,3,...,m-1), \ \  B_i=x_{i+1}\overrightarrow{P}y_i  \  \  (i=1,...,m-1),
$$
$$
l(A_i)=a_i \  \ (i=1,...,m),\  \ l(B_i)=b_i \  \ (i=1,...,m-1).
$$
Let $y\ge 0$ by an integer and
$$
a_1=a_m=\frac{y}{2}+1, \ a_2=a_3=...=a_{m-1}=0,
$$
$$
b_i=b_{m-i}=\frac{y}{2}+i+1 \ (i=1,2,...,\lfloor (m-1)/2\rfloor).
$$
If $m$ is odd, then $c=m+y+2$ and $c=\sqrt{4l+(y+1)^2}$.

If $m$ is even, we put $b_{m/2}=\frac{y}{2}+\frac{m+2}{2}$, implying that $c=m+y+2$ and $c=\sqrt{4l+(y+1)^2-1}$. Thus, the bounds in Theorem 1 are best possible.

\noindent Institute for Informatics and Automation Problems\\ National Academy of Sciences\\
P. Sevak 1, Yerevan 0014, Armenia\\ E-mail: zhora@ipia.sci.am

\end{document}